\documentclass[a4paper,11pt]{amsart}
\usepackage[arrow,matrix]{xy}
\usepackage{amsmath,amssymb,amsthm,mathrsfs,hyperref}
\theoremstyle{plain}
\newtheorem{thm}{Theorem}[section]
\newtheorem{cor}[thm]{Corollary}
\newtheorem{lem}[thm]{Lemma}
\newtheorem{prop}[thm]{Proposition}

\newtheorem{exa}[thm]{Example}

\newtheorem{rem}[thm]{Remark}

   \def\op{\oplus} \def\ot{\otimes}
\def\Hom{\operatorname {Hom}}
\def\RHom{\operatorname {RHom}}
\def\Ext{\operatorname {Ext}}

\def\End{\operatorname {End}}

\begin{document}
\title{\bf Derived $H$-module endomorphism rings}

\author{Ji-Wei He, Fred Van Oystaeyen and Yinhuo Zhang}
\address{J.-W. He\newline \indent Department of Mathematics, Shaoxing College of Arts and Sciences, Shaoxing Zhejiang 312000,
China\newline \indent Department of Mathematics and Computer
Science, University of Antwerp, Middelheimlaan 1, B-2020 Antwerp,
Belgium} \email{jwhe@usx.edu.cn}
\address{F. Van Oystaeyen\newline\indent Department of Mathematics and Computer
Science, University of Antwerp, Middelheimlaan 1, B-2020 Antwerp,
Belgium} \email{fred.vanoystaeyen@ua.ac.be}
\address{Y. Zhang\newline
\indent Department WNI, University of Hasselt,
3590 Diepenbeek, Belgium} \email{yinhuo.zhang@uhasselt.be}

\date{}

\begin{abstract} Let $H$ be a Hopf algebra, $A/B$ be an $H$-Galois extension. Let $D(A)$ and $D(B)$ be
the derived categories of right $A$-modules and of right $B$-modules
respectively. An object $M^\cdot\in D(A)$ may be regarded as an
object in $D(B)$ via the restriction functor. We discuss the
relations of the derived endomorphism rings
$E_A(M^\cdot)=\op_{i\in\mathbb{Z}}\Hom_{D(A)}(M^\cdot,M^\cdot[i])$ and
$E_B(M^\cdot)=\op_{i\in\mathbb{Z}}\Hom_{D(B)}(M^\cdot,M^\cdot[i])$. If
$H$ is a  finite dimensional semisimple Hopf algebra, then $E_A(M^\cdot)$ is a graded
subalgebra of $E_B(M^\cdot)$. In particular, if $M$ is a usual
$A$-module, a necessary and sufficient condition for $E_B(M)$ to be
an $H^*$-Galois graded extension of $E_A(M)$ is obtained. As an
application of the results, we show that the Koszul property is
preserved under Hopf Galois graded extensions.
\\[2mm]
\end{abstract}
\maketitle

\section{Introduction}

Let $A/B$ be a Hopf Galois extension over a finite dimensional Hopf
algebra $H$. The motivation of this note is to try to understand how much
homological properties, such as Koszul property and Gorenstein
property, may be preserved under Hopf Galois extensions. To this aim,
we need to discuss the derived endomorphism rings of certain
modules. The endomorphism ring extension of an $A$-module has been
studied by the second and third author in \cite{VZ}. Let $M$ be an
$A$-module. It was proved that there is an isomorphism of algebras
$\End_A(M\ot_B A)\cong\End_B(M)\#H$. Also necessary and sufficient
conditions for an endomorphism ring extension to be a Hopf Galois
extension were obtained in \cite{VZ}. We investigate whether these
properties hold when the endomorphism rings are replaced by the
derived endomorphism rings. To do this, we deal with the
relevant derived functors and work with complexes or differential graded modules instead of usual modules over an algebra $A$.

Let $H$ be a Hopf algebra with a bijective antipode, and let $A$ be
a right $H$-comodule algebra with coinvariant subalgebra
$B=A^{coH}$. We say that $A/B$ is $H$-Galois if $A$ is a Hopf Galois
extension of $B$ over $H$. Let $D(A)$ be the derived category of complexes of right $A$-modules, and $D(B)$ be the derived category
of complexes of right $B$-modules. The restriction functor
$D(A)\longrightarrow D(B)$ is an exact functor. For $M^\cdot,
N^\cdot\in D(A)$, $M^\cdot$ and $N^\cdot$ may be regarded as objects
in $D(B)$ via the restriction functor. As usual, write
$\Ext_A^*(M^\cdot,N^\cdot)$ for  $\op_{i\in\mathbb{Z}}\Hom_{D(A)}(M^\cdot,N^\cdot[i])$. Endowed with the Yoneda product
$\Ext_A^*(M^\cdot,M^\cdot)$ becomes a graded algebra. Since $A/B$ is
$H$-Galois, there is a natural $H$-action on
$\Ext^*_B(M^\cdot,M^\cdot)$ so that it becomes a graded $H$-module
algebra. In general, the graded algebra $\Ext^*_A(M^\cdot,M^\cdot)$
can not be embedded into the graded algebra
$\Ext^*_B(M^\cdot,M^\cdot)$, which differs from the usual
endomorphism rings of $A$-modules. However, if $H$ is semisimple,
then $\Ext^*_A(M^\cdot,M^\cdot)$ is the $H$-invariant subalgebra of
$\Ext^*_B(M^\cdot,M^\cdot)$. The main result of this paper is
the following one (Theorem \ref{thm3}).

\noindent{\bf Theorem.} {\it Let $H$ be a finite dimensional semisimple
Hopf algebra, and let $A/B$ be an $H$-Galois extension.
For a right $A$-module $M$, $\Ext^*_B(M,M)$ is a graded $H^*$-Galois
extension of $\Ext^*_A(M,M)$ if and only if $M\ot_BA\in \text{\rm
add}(M)$, where $\text{\rm add}(M)$ is the category consisting of
all the direct summands of finite direct sums of $M$.}

If $M$ is a right $(A,H)$-Hopf module, then $\Ext^*_B(M,M)$ is not
merely an $H^*$-Galois extension of $\Ext^*_A(M,M)$. In fact, we
have the following result (Corollary \ref{cor1}), which can be
regarded as an extension of \cite[Theorem 3.2]{Sc}.

\noindent{\bf Theorem.} {\it Let $H$ be a finite dimensional
semisimple and cosemisimple Hopf algebra. If $M$ is a right
$(A,H)$-Hopf module, then
$$\Ext^*_B(M,M)\cong\Ext^*_A(M,M)\#H^*.$$}
From this theorem, one can easily show that the Artin-Schelter Gorenstein property is preserved under Hopf Galois extensions when the Hopf algebra is semisimple and cosemisimple. With a bit more effort we can also show that the Koszul property is preserved under Hopf Galois graded extensions (Theorem \ref{thm4}).

\noindent{\bf Theorem.} {\it Let $H$ be a finite dimensional semisimple and cosemisimple Hopf algebra, let $A=\op_{n\ge0}A_n$ be a graded right $H$-module algebra such that $A_i$ is finite dimensional for all $i\ge0$, and let $B=A^{coH}$. Assume that $A/B$ is an $H$-Galois graded extension. If $B$ is an $N$-Koszul algebra, then $A$ is an $N$-Koszul algebra.}

In particular, for a pointed cocommutative Hopf algebra $A$, the
graded algebra $gr(A)$ associated to the coradical filtration is a
Koszul algebra if $A$ has finitely many group-like elements and the
dimension of $P(A)$, the vector space of primitive elements, is
finite.

The main tools used in this note are differential graded algebras
and differential graded modules. Let us recall some notations and properties of differential graded (dg, for short) algebras and dg modules which will
be used later. By a dg algebra we mean a graded algebra
$R=\op_{n\in\mathbb{Z}}R^n$ with a differential $d$ of degree 1 which is
also a derivation of $R$. Similarly we have dg modules. We sometimes
view a usual algebra $A$ as a dg algebra concentrated in degree 0
with a trivial differential. Then a complex of $A$-modules can be
regarded as a dg $A$-module. Let $R$ be a dg algebra, and let $X_R$
and $Y_R$ be right dg $R$-modules. Write $\Hom_R(X,Y)$ for $\op_{i\in\mathbb{Z}}\Hom_R^i(X,Y)$, where $\Hom_R^i(X,Y)$ is the set of all  graded $R$-module morphisms of degree $i$. $\Hom_R(X,Y)$ is a complex of
vector space. The differential on $\Hom_R(X,Y)$ is given by
$\delta(f)=d_Y\circ f-(-1)^{|f|}f\circ d_X$, where $|f|$ is the
degree of $f$. Also $\Hom_R(X,X)$ is a dg algebra and $\Hom_R(X,Y)$
is a right dg $\Hom_R(X,X)$-module. Let $H$ be a Hopf algebra. We
say that $R$ is a {\it dg left $H$-module algebra} if $R$ is a left
$H$-module algebra and the $H$-action respects the the grading of
$R$ and is compatible with the differential of $R$. In this case,
the smash product $R\#H$ is a dg algebra and the cohomology algebra
$H^\cdot(R\#H)\cong H^\cdot(R)\#H$. Notation and terminology used in
this note is the same as in \cite{AFH,FHT,KZ}. For example, we use $\RHom_R(-,-)$ to denote the right derived functor of $\Hom_R(-,-)$, and $-\ot^L_R-$ to denote the left derived functor of the tensor functor $-\ot_R-$ over a dg algebra $R$.

We work over a fixed field $k$. All the algebras and modules are
defined over $k$. Unadorned $\ot$ means $\ot_k$, and $\Hom$ means
$\Hom_k$.

\section{Derived endomorphism rings}

Throughout $H$ will be always a Hopf algebra with a bijective
antipode. Let $A/B$ be an $H$-Galois extension with the canonical
map $\beta:A\ot_BA\longrightarrow A\ot H$. For any $h\in H$, we write
$\sum X_i^h\ot Y_i^h$ for $\beta^{-1}(1\ot h)$. For right $A$-modules $M$
and $N$, there is a natural left $H$-module action on $\Hom_B(M,N)$:

For $h\in H$, $m\in M$ and $f\in\Hom_B(M,N)$, $$(h\cdot f)(m)=\sum
f(mX_i^{Sh})Y_i^{Sh}.$$ Similarly, if $M^\cdot$ and $N^\cdot$ are
complexes of right $A$-modules, then $\Hom_B(M^\cdot,N^\cdot)$ is a
complex of left $H$-modules. Hence for any $M^\cdot\in D^-(A)$ and
$N^\cdot\in D(A)$, $\RHom_B(M^\cdot,N^\cdot)$ is a complex of left
$H$-modules.

If $M$ is a right $A$-module, then there is an isomorphism of right
$B$-modules and right $H$-comodules:
$$M\ot_BA\longrightarrow M\ot H,\qquad m\ot a\mapsto\sum ma_{0}\ot
a_{1},$$ with the inverse given by:
$$M\ot H\longrightarrow M\ot_B
A,\qquad m\ot h\mapsto\sum mX_i^h\ot Y_i^h.$$
These isomorphisms may
be generalized to complexes. If $M^\cdot$ is a complex of
$A$-modules, then we have an isomorphism of complexes:
$M^\cdot\ot_BA\longrightarrow M^\cdot\ot H$.

\begin{thm} \label{thm1}Let $A/B$ be an $H$-Galois extension, and $M^\cdot$ an object in $D^-(A)$. If $H$ is finite dimensional, then there is
an isomorphism of graded algebras
$$\Ext^*_A(M^\cdot\ot_B^LA,M^\cdot\ot_B^LA)\longrightarrow
\Ext^*_B(M^\cdot,M^\cdot)\#H.$$
\end{thm}
\proof Let $P^\cdot$ be a bounded above complex of projective
$A$-modules which is quasi-isomorphic to the complex $M^\cdot$. Then
$M^\cdot\ot_B^LA\cong P^\cdot\ot_B A$. Of course $P^\cdot\ot_BA$ is
homotopically projective (see \cite[Ch. 8]{KZ}). Hence
$$
\begin{array}{ccl}
  \RHom_A(M^\cdot\ot_B^LA,M^\cdot\ot_B^LA) & \cong & \RHom_A(P^\cdot\ot_BA,P^\cdot\ot_BA) \\
   & = &\Hom_A(P^\cdot\ot_BA,P^\cdot\ot_BA)\\
   & \cong&\Hom_B(P^\cdot,P^\cdot\ot_BA)\\
   &\cong&\Hom_B(P^\cdot,P^\cdot\ot H)\\
   &\cong&\Hom_B(P^\cdot,P^\cdot)\ot H\\
   &=&\RHom_B(M^\cdot,M^\cdot)\ot H.
\end{array}
$$
Therefore as graded spaces, we have
$\Ext^*_A(M^\cdot\ot_B^LA,M^\cdot\ot^L_BA)\cong\Ext^*_B(M^\cdot,M^\cdot)\#
H$. It remains to show that this isomorphism is an isomorphism
of graded algebras.

From the above isomorphisms, we get an isomorphism of complexes
\begin{eqnarray}\label{eq1}
\varphi:\Hom_B(P^\cdot,P^\cdot)\ot
H&\longrightarrow&\Hom_A(P^\cdot\ot_BA,P^\cdot\ot_BA),
\end{eqnarray} such that for
$f\in\Hom_B(P^\cdot,P^\cdot)$, $h\in H$, $p\in P$ and $a\in A$,
$$\varphi(f\ot h)(p\ot_B a)=\sum f(p)X_i^h\ot_B Y_i^ha.$$ We have to
show that $\varphi$ is an isomorphism of dg algebras from
$\Hom_B(P^\cdot,P^\cdot)\# H$ to
$\Hom_A(P^\cdot\ot_BA,P^\cdot\ot_BA)$. Since $\varphi$ is compatible
with the differentials, it suffices to show that it is an algebra
morphism. The proof is essentially the same as that of \cite[Theorem
2.3]{VZ}, so we omit it. Now by taking the cohomology algebra of the
dg algebras $\Hom_B(P^\cdot,P^\cdot)\# H$ and
$\Hom_A(P^\cdot\ot_BA,P^\cdot\ot_BA)$ we get the desired result.
\qed

Let $M^\cdot$ be a complex of $A$-modules. We have
$\Hom_A(M^\cdot,M^\cdot)\cong\Hom_B(M^\cdot,M^\cdot)^H$ as dg
algebras, where $\Hom_B(M^\cdot,M^\cdot)^H$ is the $H$-invariant dg
subalgebra of $\Hom_B(M^\cdot,M^\cdot)$. If $H$ is a  finite dimensional semisimple Hopf algebra, then it is not difficult to see
that $\Ext_A^*(M^\cdot,M^\cdot)\cong\Ext^*_B(M^\cdot,M^\cdot)^H$ as graded algebras. In
other words, the graded algebra $\Ext_B^*(M^\cdot,M^\cdot)$ is an
$H^*$-extension of $\Ext_A^*(M^\cdot,M^\cdot)$. We ask when this
$H^*$-extension is an $H^*$-Galois extension.

Let $M^\cdot$ be a complex of right $(A,H)$-Hopf modules. Then
$\Hom_A(M^\cdot,M^\cdot)$ is a dg left $H^*$-module algebra with the
left $H^*$-action given by, for $\alpha\in H^*$ and
$f\in\Hom_A(M^\cdot,M^\cdot)$,
$$(\alpha\cdot f)(m)=\sum \alpha_{(1)}f(S(\alpha_{(2)})m),$$ which induces a left $H^*$-module algebra structure on
$\Ext_A^*(M^\cdot,M^\cdot)$.

\begin{thm}\label{thm2} Let $H$ be a   finite dimensional semisimple and cosemisimple Hopf algebra.
Let $A/B$ be an $H$-Galois extension. If $M^\cdot$ is a bounded
 above complex of right $(A,H)$-Hopf modules, then there is an isomorphism
of graded algebras
$$\Ext^*_B(M^\cdot,M^\cdot)\cong\Ext^*_A(M^\cdot,M^\cdot)\#H^*.$$
\end{thm}
\proof The complex $M^\cdot$ of right $(A,H)$-Hopf modules may be
regarded as a complex of right $A\#H^*$-modules in a natural way.
By \cite[Theorem 2.2]{CFM} the functor
$\Hom_B(A,-):\text{Mod-}B\longrightarrow \text{Mod-}A\#H^*$ is an
equivalence of abelian categories with inverse functor
$-\ot_{A\#H^*}A$. Hence we have
$$\begin{array}{cl}
    &\Ext^*_{A\#H^*}(M^\cdot\ot H^*,M^\cdot\ot H^*)    \\
     \cong & \Ext_B^*\left((M^\cdot\ot
H^*)\ot_{A\#H^*}A,(M^\cdot\ot
H^*)\ot_{A\#H^*}A\right) \\
    \cong &\Ext_B^*(M^\cdot,M^\cdot).
  \end{array}
$$ The last isomorphism holds because $(M^\cdot\ot
H^*)\ot_{A\#H^*}A\cong(M^\cdot\ot_A(A\#H^*))\ot_{A\#H^*}A
\cong M^\cdot$ as a complex of right $B$-modules. On the other hand, since
as complexes of right $A\#H^*$-modules
$$M^\cdot\ot H^*\cong M^\cdot\ot_A^L(A\#H^*),$$ we have
$$\begin{array}{cl}
 & \Ext^*_{A\#H^*}(M^\cdot\ot H^*,M^\cdot\ot H^*) \\
 \cong & \Ext^*_{A\#H^*}(M^\cdot\ot_A^L (A\#H^*),M^\cdot\ot_A^L
(A\#H^*)) \\
\cong & \Ext^*_{A}(M^\cdot,M^\cdot)\#H^*.
\end{array}
$$
The last isomorphism follows from Theorem \ref{thm1} since $A\#H^*$
is an $H^*$-Galois extension of $A$. Now we obtain the desired
isomorphism
$$\Ext_B^*(M^\cdot,M^\cdot)\cong \Ext^*_{A}(M^\cdot,M^\cdot)\#H^*.
\qed$$

From the above theorem, we obtain the following corollary which may
be regarded as a natural generalization of \cite[Theorem 3.2]{Sc}.

\begin{cor}\label{cor1} Let $H$ be a   finite dimensional semisimple and cosemisimple Hopf algebra, and let
$M$ be a right $(A,H)$-Hopf module. Then
$$\Ext^*_B(M,M)\cong\Ext^*_A(M,M)\#H^*.$$
\end{cor}

For an $A$-module $M$, the algebra $\Ext_B^*(M^\cdot,M^\cdot)$ may
not be a smash product of $\Ext_A^*(M^\cdot,M^\cdot)$ with $H^*$.
However, if the Hopf algebra $H$ is unimodular, then we will see
that in certain cases the graded algebra $\Ext_B^*(M^\cdot,M^\cdot)$
is an $H^*$-Galois graded extension of $\Ext_A^*(M^\cdot,M^\cdot)$.
Firstly, we make explicit the exact meaning of a Galois graded
extension.

Let $E$ be a $\mathbb{Z}$-graded algebra, $H$ a Hopf algebra. We say
that $E$ is a graded right $H$-comodule algebra if $E$ is a right
$H$-comodule algebra and the $H$-coaction on $E$ respects the
grading of $E$, that is; for $x\in E_n$, $\rho(x)\in E_n\ot H$. Let
$D=E^{coH}$. Then $D$ is a graded subalgebra of $E$. $E/D$ is called
an {\it $H$-Galois graded extension} if the canonical map \cite{Sch}
$$E\ot_D E\longrightarrow E\ot H,\quad x\ot x'\mapsto \sum
xx'_{(0)}\ot x'_{(1)}$$ is bijective.

We have the following observation.

\begin{lem} \label{lem1} Let $E/D$ be an $H$-Galois graded extension. If $E=\op_{n\ge0} E_n$ is positively
graded, then $E_0/D_0$ is an $H$-Galois extension.
\end{lem}

For graded right $E$-modules $V$ and $W$, let
$\underline{\Hom}_E(V,W)=\op_{i\in\mathbb{Z}}\text{Hom}_E^i(V,W)$, where
$\text{Hom}_E^i(V,W)$ is the set of all the graded $E$-module
morphisms of degree $i$, and $\underline{\text{End}}(V_E)$ denotes
the sum of all the graded endomorphisms. Then
$\underline{\Hom}_E(V,W)$ is a graded vector space and
$\underline{\text{End}}(V_E)$ is a graded algebra. We use $\underline{\Ext}_E^*(-,-)$ to denote the derived functor of $\underline{\Hom}_E(-,-)$. Note that $\underline{\Ext}^n_E(V,W)$ is a graded vector space for all $n\in \Bbb Z$, where the grading is induced from the gradings of $V$ and $W$. Hence $\underline{\Ext}^*_E(V,W)$ is a bigraded vector space. To avoid possible confusions, following \cite{GMMZ}, we say an element $x$ in $\underline{\Ext}^n_E(V,W)$ has {\it ext-degree} $n$ when we ignore the grading induced from the gradings of $V$ and $W$.

For $H$-Galois graded extensions, we have the following result which
corresponds to \cite[Theorem 1.2]{CFM}. The proof is exactly
the same as that of \cite[Theorem 1.2]{CFM} except that we should
keep in mind that everything needs to preserve the gradings.

\begin{lem}\label{lem4} Let $H$ be a finite dimensional Hopf
algebra. The following are equivalent:
\begin{itemize}
  \item [(i)] $E/D$ is an $H$-Galois graded extension.
  \item [(ii)] {\rm(a)} The natural map $E\#H^*\longrightarrow
  \underline{\text{\rm End}}(E_D)$ is an isomorphism of graded
  algebras.\\
  {\rm(b)} $E_D$ is a finitely generated projective graded $D$-module.
  \item [(iii)] As a left graded $E\#H^*$-module, $E$ is a generator
  of the category of graded $E\#H^*$-modules.
\end{itemize}
\end{lem}

\begin{rem}\label{rem1} {\rm If the $H$-Galois extension $A/B$ in Theorem \ref{thm1} (resp. \ref{thm2}) is $H$-Galois graded extension and $M^\cdot$ is a bounded above complex of graded $A$-modules (resp. of graded right $(A,H)$-Hopf modules), then the result still holds when the functor $\Ext$ is replaced by $\underline{\Ext}$. Moreover, one can check that the isomorphisms in the theorems are bigraded morphisms. In particular, the isomorphism in Corollary \ref{cor1} has the following form, which we will use later, $$\underline{\Ext}^*_B(M,M)\cong\underline{\Ext}^*_A(M,M)\#H^*,$$ where the bigrading of the right hand is induced from the bigrading of $\underline{\Ext}^*_A(M,M)$.}
\end{rem}

For later use, we need to generalize some results in \cite{D,VZ}.

\begin{lem}\label{lem2} Let $H$ be a finite dimensional Hopf algebra, and let $A/B$ be an $H$-Galois
extension. If $M^\cdot\in D^-(A)$ and $N^\cdot$ is a complex of
right $B$-modules, then there is an isomorphism in $D(k)$
$$\RHom_B(M^\cdot,N^\cdot)\cong\RHom_A(M^\cdot,N^\cdot\ot^L_BA).$$
\end{lem}
\proof Replace $M^\cdot$ and $N^\cdot$  by complexes $P^\cdot$ and $Q^\cdot$ of projective $A$-modules respectively. Since $A/B$ is an $H$-Galois
extension, $A_B$ is a projective module. Hence $P^\cdot$ and
$Q^\cdot$ are complexes of projective $B$-modules. We have
$\RHom_B(M^\cdot,N^\cdot)=\Hom_B(P^\cdot,Q^\cdot)$, and
$\RHom_A(M^\cdot,N^\cdot\ot^L_BA)=\Hom_A(P^\cdot,Q^\cdot\ot_BA)$.\\
Following \cite{D}, we define a map
\begin{equation}\label{eq2}
    \xi:\Hom_B(P^\cdot,Q^\cdot)\longrightarrow\Hom_A(P^\cdot,Q^\cdot\ot_BA)
\end{equation}
by $\xi(f)(p)=\sum f(pX^t_i)\ot_B Y^t_i$ where $t$ is a nonzero
right integral in $H$. We claim that $\xi$ is an isomorphism of
complex of vector spaces. By essentially the same calculations as in
the proof of \cite[Theorem 5]{D}, one sees that $\xi$ is bijective.
What left to show is that $\xi$ is compatible with the
differentials. Now for $f\in \Hom_B(P^\cdot,Q^\cdot)$ and $p\in
P^\cdot$, we have
$$\begin{array}{ccl}
    \xi(d(f))(p) & = & \xi(d_{Q^\cdot}\circ f-(-1)^{|f|}f\circ d_{P^\cdot})(p) \\
     & = & \sum d_{Q^\cdot}(f(pX_i^t))\ot_B
Y_i^t-\sum(-1)^{|f|}f(d_{P^\cdot}(pX_i^t))\ot_B Y_i^t \\
    & = & d_{Q^\cdot\ot_B A}(\sum f(pX_i^t)\ot_B
Y_i^t)-(-1)^{|f|}\sum f(d_{P^\cdot}(p)X_i^t)\ot_B Y_i^t\\
&=&d_{Q^\cdot\ot_B A}\circ \xi(f)(p)-(-1)^{|f|} \xi(f)\circ
d_{P^\cdot}(p)\\
&=&d(\xi(f))(p).
  \end{array}
$$
Hence we get the desired result. \qed

Let $H$ be a   finite dimensional semisimple Hopf algebra. By Theorem
\ref{thm1}, we have $K=\Ext^*_A(M^\cdot\ot^L_B A,M^\cdot\ot^L_B
A)\cong\Ext^*_B(M^\cdot,M^\cdot)\#H$. Write
$E=\Ext^*_B(M^\cdot,M^\cdot)$ and $D=\Ext^*_A(M^\cdot,M^\cdot)$.
Then $E$ is a graded left $K$-module and a graded right $D$-module.
Since $\Hom_B(P^\cdot,P^\cdot)$ is a dg $H$-module algebra, it is a
left dg $\Hom_B(P^\cdot,P^\cdot)\#H$-module. Let
$U=\Hom_B(P^\cdot,P^\cdot)^H=\Hom_A(P^\cdot,P^\cdot)$ be the fixed
dg subalgebra of $\Hom_B(P^\cdot,P^\cdot)$. Then
$\Hom_B(P^\cdot,P^\cdot)$ is a dg
$\Hom_B(P^\cdot,P^\cdot)\#H$-$U$-bimodule. By taking the cohomology,
we obtain that $E=H^\cdot\Hom_B(P^\cdot,P^\cdot)$ is a graded left
$K$- and right $D$-bimodule since $K\cong E\#H$ and $D\cong H^\cdot
U$.

\begin{prop}\label{prop1} With notation as above, there is a graded bimodule
isomorphism
$${}_KE_D\cong {}_K\Ext^*_B(M^\cdot,M^\cdot\ot^L_BA)_D.$$
\end{prop}
\proof By Lemma \ref{lem2},
$E\cong\Ext^*_B(M^\cdot,M^\cdot\ot^L_BA)$ as graded spaces. It
suffices to prove that the isomorphism is a bimodule isomorphism.
Going back to the morphism defined by (\ref{eq2}) in the proof of
Lemma \ref{lem2}, one sees that $\xi$ is a dg right
$\Hom_A(P^\cdot,P^\cdot)$-module morphism. We show that
$\xi:\Hom_B(P^\cdot,P^\cdot)\longrightarrow\Hom_A(P^\cdot,P^\cdot\ot_BA)$
is a left dg $\Hom_B(P^\cdot,P^\cdot)\#H$ module morphism. Note that
the left dg $\Hom_B(P^\cdot,P^\cdot)\#H$ module action is given by
the pullback of dg algebra morphism $\varphi$ defined by (\ref{eq1})
in the proof of Theorem \ref{thm1}. Now for
$g,f\in\Hom_B(P^\cdot,P^\cdot)$, $h\in H$ and $p\in P^\cdot$, we
have $$\begin{array}{ccl}
         \xi((g\#h)\cdot f)(p) & = & \xi(g\circ(h\cdot f))(p) \\
          & = & \sum g\circ(h\cdot f)(pX_i^t)\ot_B Y^t_i \\
         & = & \sum g(f(pX_i^tX^{Sh}_j)Y_j^{Sh})\ot_B Y_i^t.
       \end{array}
$$
Since $A/B$ is an $H$-Galois extension and $H$ is semisimple, it is
not hard to check (see \cite{VZ}) $\sum X_i^tX_j^{Sh}\ot_B
Y_j^{Sh}\ot_B Y_i^t=\sum X_i^t\ot_B X_j^{h}\ot_B Y^h_jY^t_i$ in
$A\ot_BA\ot_BA$. Then we obtain
$$\begin{array}{ccl}
     \xi((g\#h)\cdot f)(p) & = & \sum g(f(pX_i^t)X^{h}_j))\ot_B Y^h_jY_i^t \\
     & = & \varphi(g\#h)\circ\xi(f)(p) \\
     & = & (g\#h)\cdot\xi(f)(p).
  \end{array}
$$
Therefore $\xi$ is a dg bimodule isomorphism. By taking the
cohomologies of $\Hom_B(P^\cdot,P^\cdot)$ and
$\Hom_A(P^\cdot,P^\cdot\ot_B A)$ respectively we arrive at the
desired graded bimodule isomorphism. \qed

Let $\mathcal{T}$ and $\mathcal{D}$ be triangulated categories. Let
$F:\mathcal{T}\longrightarrow\mathcal{D}$ and
$G:\mathcal{D}\longrightarrow\mathcal{T}$ be a pair of adjoint exact
functors with $F$ left adjoint to $G$. Recall that the Auslander
class is the subcategory of $\mathcal{T}$:
$$\mathcal{A}(\mathcal{T})=\{X\in \mathcal{T}|\text{the adjunction
map}\ X\longrightarrow GFX\ \text{is isomorphic}\},$$ and the Bass
class is the subcategory of $\mathcal{D}$:
$$\mathcal{B}(\mathcal{D})=\{Y\in \mathcal{D}|\text{the adjunction map}\ FGY\longrightarrow Y\ \text{is
isomorphic}\}.$$ It is well known that both the Auslander class and
Bass class are thick triangulated subcategories of the respective
triangulated categories, and the adjoint functors $F$ and $G$ induce
a pair of inverse equivalences between triangulated categories
$\mathcal{A}(\mathcal{T})$ and $\mathcal{B}(\mathcal{D})$.

The following lemma is straightforward.
\begin{lem}\label{lem3} For $M^\cdot\in D^-(A)$, if $M^\cdot\ot^L_BA$ is a direct summand of finite direct sum of $M^\cdot$,
then $\Ext_A^*(M^\cdot,M^\cdot\ot^L_BA)$ is a finitely generated
graded projective $\Ext^*_A(M^\cdot,M^\cdot)$-module.
\end{lem}

Let $P^\cdot\in D^-(A)$ be a complex of projective $A$-modules. We
regard $A$ as a dg algebra concentrated in degree 0. Then $P^\cdot$
is a dg right $A$-module. Let $R=\Hom_A(P^\cdot,P^\cdot)$ be the dg
algebra of endomorphisms. Then $P^\cdot$ is a left dg $R$-module,
and it is also a dg $R$-$A$-bimodule. Therefore we have an exact
functor
$$\RHom_A(P^\cdot,-):D(A)\longrightarrow D_{dg}(R),$$ where
$D_{dg}(R)$ is the derived category of right dg $R$-modules. This
functor is naturally left adjoint to the functor
$$-\ot_R^LP^\cdot:D_{dg}(R)\longrightarrow D(A).$$ Clearly the dg
$A$-module $P^\cdot$ lies in the Auslander class. Since the
Auslander class is a thick triangulated subcategory, all the direct
summands of finite direct sums of copies of $P^\cdot$ belong to the
Auslander class.

Let $M$ be a right $A$-module. Denote by $\text{add}(M)$
the category of $A$-modules isomorphic to direct summands of
finite sums of $M$. Recall that $E=\Ext_B^*(M,M)$, $D=\Ext_A^*(M,M)$ and
$K=E\#H$.

\begin{thm}\label{thm3} Let $H$ be a  finite dimensional semisimple  Hopf algebra, and
$A/B$ be an $H$-Galois extension. For a right $A$-module $M$, the
following are equivalent.

{\rm(i)} $M\ot_BA\in \text{\rm add}(M)$.

{\rm(ii)} $\End_B(M)/\End_A(M)$ is an $H^*$-Galois extension.

{\rm(iii)} $E/D$ is an $H^*$-Galois graded extension.

Moreover, if $M\ot_BA\in \text{\rm add}(M)$, then $D$ and $E\#H$ are
graded Morita equivalent.
\end{thm}
\proof The equivalence between (i) and (ii) was proved in \cite[Theorem
2.4]{VZ}. It remains to show that (i) and (iii) are equivalent.

Assume that $E/D$ is an $H^*$-Galois graded extension. Since
$E=\Ext_B^*(M,M)$ is positively graded, by Lemma \ref{lem1},
$\Ext_B^0(M,M)$ is an $H^*$-Galois extension of $\Ext^0_B(M,M)$,
that is; End$_B(M)$/End$_A(M)$ is an $H^*$-Galois extension. Hence
$M\ot_BA\in \text{add}(M)$.

Conversely, assume $M\ot_BA\in\text{add}(M)$. By Proposition
\ref{prop1}, ${}_KE_D\cong {}_K\Ext^*_B(M,M\ot^L_BA)_D$ as graded
$K$-$D$-bimodules. It follows from Lemma \ref{lem3} that $\Ext^*_B(M,M\ot^L_BA)$, and
hence $E$, is a finitely generated graded right $D$-module.

Now we prove that $\underline{\text{End}}(E_D)$ is isomorphic to the
graded algebra $K=E\#H$. Let $P^\cdot\in D^-(A)$ be a projective
resolution of $M_A$. Viewing $A$ as a dg algebra concentrated in
degree 0, we see that $P^\cdot$ is a dg $R$-$A$-bimodule, where
$R=\Hom_A(P^\cdot,P^\cdot)$. As mentioned before,
$(\RHom_A(P^\cdot,-),-\ot_R^LP^\cdot)$ is a pair of adjoint
functors; and $M\cong P^\cdot$ lies in the Aulander class of the
adjoint pair. Therefore $\text{add}(M)$ is contained in the
Auslander class. We have the following graded algebra isomorphism:
$$\begin{array}{cl}
     & \bigoplus_{i\in\mathbb{Z}}\Hom_{D_{dg}(R)}(\RHom_A(P^\cdot,M\ot_BA),
     \RHom_A(P^\cdot,M\ot_BA)[i]) \\
    \cong &\bigoplus_{i\in\mathbb{Z}} \Hom_{D(A)}(M\ot_BA,M\ot_BA[i]).
  \end{array}
$$
Since $P^\cdot$ is a projective resolution,
$\RHom_A(P^\cdot,M\ot_BA)\cong\Hom_A(P^\cdot,P^\cdot\ot_BA)$. Since
$M\ot_BA\in\text{add}(M)$, $P^\cdot\ot_BA$ is homotopically
equivalent to $Q^\cdot$. Thus there exist an integer $n$ and some complex $\bar{Q^\cdot}$ such  that $Q^\cdot\op
\bar{Q^\cdot}$ is homotopically equivalent to $(P^\cdot)^{\op n}$. It follows that $\Hom_A(P^\cdot,P^\cdot\ot_BA)\cong\Hom_A(P^\cdot,Q^\cdot)$.
Moreover, $\Hom_A(P^\cdot,Q^\cdot)$ is a homotopically projective dg
$R$-module. Now we have the following graded algebra isomorphisms:
$$\begin{array}{cl}
     & \bigoplus_{i\in\mathbb{Z}}\Hom_{D_{dg}(R)}(\RHom_A(P^\cdot,M\ot_BA),\RHom_A(P^\cdot,M\ot_BA)[i]) \\
     \cong &H^\cdot\RHom_{R}(\RHom_A(P^\cdot,M\ot_BA),\RHom_A(P^\cdot,M\ot_BA)) \\
   \cong &H^\cdot
   \Hom_R(\Hom_A(P^\cdot,Q^\cdot),\Hom_A(P^\cdot,Q^\cdot))\\
   \overset{(a)}\cong&\underline{\Hom}_{H^\cdot
   R}(H^\cdot\Hom_A(P^\cdot,Q^\cdot),H^\cdot\Hom_A(P^\cdot,Q^\cdot))\\
  =&\underline{\Hom}_{D}(\Ext^*_A(M,M\ot_BA),\Ext^*_A(M,M\ot_BA))\\
  \cong&\underline{\text{End}}(E_D).
  \end{array}
$$ The isomorphism $(a)$ holds because $\Hom_A(P^\cdot,Q^\cdot)$ is
a direct summand (in the homotopy category) of a free dg $R$-module.
More precisely, in the homotopy category of right dg $R$-modules, we
have
$\Hom_A(P^\cdot,Q^\cdot)\op\Hom_A(P^\cdot,\bar{Q}^\cdot)\cong\Hom_A(P^\cdot,Q^\cdot\op
\bar{Q}^\cdot)\cong\Hom_A(P^\cdot,{P^\cdot}^{\op n})\cong R^{\op
n}$. Let $X=\Hom_A(P^\cdot,Q^\cdot)$ and
$Y=\Hom_A(P^\cdot,\bar{Q}^\cdot)$. We have the following commutative
diagram: $$\xymatrix{
  H^\cdot\Hom_R(X\op Y,X) \ar[d]_{\cong} \ar[r] & \Hom_{H^\cdot R}(H^\cdot(X\op Y),H^\cdot X) \ar[d]^{\cong} \\
    H^\cdot\Hom_R(R^{\op n},X)\ar[r]^{\cong\quad} &  \Hom_{H^\cdot R}(H^\cdot(R^{\op n}),H^\cdot X).  }$$
Since all the morphisms in the diagram are natural, we obtain a natural  isomorphism $H^\cdot\Hom_R(X\op
Y,X)\cong\Hom_{H^\cdot R}(H^\cdot(X\op Y),H^\cdot X)$, which implies
the isomorphism $(a)$: $H^\cdot\Hom_R(X,X)\cong\Hom_{H^\cdot
R}(H^\cdot X,H^\cdot X)$.

On the other hand, we have $$\bigoplus_{i\in\mathbb{Z}}
\Hom_{D(A)}(M\ot_BA,M\ot_BA[i])\cong\Ext^*_A(M\ot_BA,M\ot_BA)\cong
E\#H.$$ Thus $\underline{\text{End}}(E_D)\cong E\#H$  as graded
algebras. Now applying Lemma \ref{lem4}, we obtain that $E/D$ is an
$H^*$-Galois extension.

Assume now $M\ot_BA\in \text{add}(M)$. Then $E/D$ is an $H^*$-Galois
extension. By Lemma \ref{lem4}(iii), ${}_{E\#H}E$ is a graded
generator. Since $H$ is semisimple, $E$ is a graded
projective ${E\#H}$-module. Therefore, $E$ is finitely generated
projective generator of the category of graded ${E\#H}$-modules.
Thus $E\#H$ is graded Morita equivalent to $D$ (see \cite{GG}). \qed

We next want to show that Galois extensions over a semisimple and cosemisimple Hopf algebra preserve Koszul property. Let us recall the definition of Koszul algebras. Let $A=\op_{n\ge0}A_n$ be a graded algebra such that $A_i$ is finite dimensional for all $i\ge0$, $A_0$ is semisimple and $A_iA_j=A_{i+j}$ for all $i,j\ge0$. Let $N\ge2$ be an integer. A graded $A$-module $M$ is called an {\it $N$-Koszul module} (see, \cite{GMMZ} for instance) if $M$ has a graded projective resolution: $$\cdots\longrightarrow P^{-n}\longrightarrow\cdots\longrightarrow P^{-1}\longrightarrow P^0\longrightarrow M\longrightarrow0$$
such that the graded projective module $P^{-n}$ ($n>0$) is generated in degree $\frac{n}{2}N$ if $n$ is even, or $\frac{n-1}{2}N+1$ if $n$ is odd. If the trivial graded right $A$-module $A_0$ (via the projection) is $N$-Koszul, then $A$ is called an {\it $N$-Koszul algebra}.

\begin{thm}\label{thm4} Let $H$ be a finite dimensional semisimple and cosemisimple Hopf algebra, let $A=\op_{n\ge0}A_n$ be a graded right $H$-module algebra such that $A_i$ is finite dimensional for all $i\ge0$, and let $B=A^{coH}$. Assume that $A/B$ is an $H$-Galois graded extension. If $B$ is an $N$-Koszul algebra, then $A$ is an $N$-Koszul algebra.
\end{thm}
\proof By the assumption, $B_0$ is a finite dimensional semisimple algebra. Since $A/B$ is an $H$-Galois graded extension, $A_0/B_0$ is an $H$-Galois extension by Lemma 2.4. Now that $H$ is semisimple implies that $A_0\# H$ is semisimple since it is Morita equivalent with $B_0$. Since $H$ is cosemisimple, $(A_0\# H)\# H^*$ is semisimple by Maschke's theorem. It follows from the Morita equivalence between $A_0$ and $(A_0\# H)\# H^*$ that $A_0$ is a semisimple algebra. As a right $B_0$-module, $A_0=B_0\op S$ for some finite dimensional $B_0$-module $S$. By Remark \ref{rem1}, we have an isomorphism of bigraded algebras
\begin{equation}\label{eq3}
    \underline{\Ext}^*_B(A_0,A_0)\cong\underline{\Ext}^*_A(A_0,A_0)\#H.
\end{equation}
Since $B$ is $N$-Koszul, the graded space $\underline{\Ext}^1_B(B_0,B_0)$ is concentrated in degree 1, and $\underline{\Ext}^2_B(B_0,B_0)$ is concentrated in degree $N$ (see \cite{GMMZ}). Therefore $\underline{\Ext}^1_B(B_0,S)$, ($\underline{\Ext}^1_B(S,B_0)$ and $\underline{\Ext}^1_B(S,S)$, respectively) is concentrated in degree 1, and $\underline{\Ext}^2_B(B_0,S)$, ($\underline{\Ext}^2_B(S,B_0)$ and $\underline{\Ext}^2_B(S,S)$ respectively) is concentrated in degree $N$ as $S$ is a direct summand of of a finite sum of $B_0$. Hence $\underline{\Ext}^1_B(A_0,A_0)\cong\underline{\Ext}^1_B(B_0,B_0)\op\underline{\Ext}^1_B(B_0,S)\op\underline{\Ext}^1_B(B_0,S)
\op\underline{\Ext}^1_B(S,S)$ is concentrated in degree 1. Since we already know that $A_0$ is semisimple, $A$ must be generated in degrees 0 and 1. Similarly, we see that $\underline{\Ext}^2_B(A_0,A_0)$ is concentrated in degree $N$. Thus $A$ must be homogeneous in the sense of \cite{GMMZ}.

By the isomorphism (\ref{eq3}), the graded algebra $\underline{\Ext}^*_A(A_0,A_0)$ is generated in ext-degrees 0, 1 and 2 if and only if $\underline{\Ext}^*_B(A_0,A_0)$ is. Next we show that $\underline{\Ext}^*_B(A_0,A_0)$ is generated in ext-degrees 0, 1 and 2. For convenience, we let $E^n=\underline{\Ext}^n_B(B_0,B_0)$, $D^n=\underline{\Ext}^n_B(S,S)$, $U^n=\underline{\Ext}^n_B(S,B_0)$ and $V^n=\underline{\Ext}^n_B(B_0,S)$. Write $E$ and $D$ for the graded algebras $\op_{n\ge0}E^n$ and $\op_{n\ge0}D^n$ respectively. Similarly let  $U$ and $V$ be the respective graded $E$-$D$-bimodule $\op_{n\ge0}U^n$ and $D$-$E$-bimodule $\op_{n\ge0}V^n$.  Since $S$ as a graded left $E$-module is a direct summand of of a finite sum of $B_0$, $U$ is generated in ext-degree 0. Similarly, $V$ is generated in ext-degree 0 as a graded right $E$-module. The graded algebra $\underline{\Ext}^*_B(A_0,A_0)$ is isomorphic to the following matrix algebra $$\left(
                                                                                  \begin{array}{cc}
                                                                                    E & U \\
                                                                                    V & D \\
                                                                                  \end{array}
                                                                                \right)
$$ where the products $UV$ and $VU$ are the Yoneda products of the extensions. Now if $B$ is $N$-Koszul, then $E$ is generated by $E^0$, $E^1$ and $E^2$. Since $V$, as a graded right $E$-module, is generated in degree 0, we have $V=V^0E$. Hence each element of $V$ can be written as a sum of some multiplications of elements in $E^0$, $E^1$, $E^2$ and $V^0$. Similarly, each element of $U$ can be written as a sum of some multiplications of elements in $E^0$, $E^1$, $E^2$ and $U^0$. Since $B_0$ is semisimple, we may assume that there is a semisimple $B_0$-module $T$ such that $S\op T\cong \op_{i=1}^n Q_i$ as right $B_0$-modules, where $Q_i=B_0$ for all $i=1,\dots, n$. Notice that the actions of an element in $\op_{i\ge1}B_i$ on both sides of the foregoing isomorphism are trivial. Thus the isomorphism is in fact a right $B$-module isomorphism. Let $\iota:S\longrightarrow \op_{i=1}^n Q_i$ be the inclusion map, and $\pi:\op_{i=1}^n Q_i\longrightarrow S$ be the projection map. Then $\pi\circ\iota=id$. For any $n\ge0$, let $\iota^n_{ext}:\underline{\Ext}^n_B(S,S)\longrightarrow \underline{\Ext}^n_B(S,\op_{i=1}^nQ_i)$ be the map induced by $\iota$, and $\pi^n_{ext}:\underline{\Ext}^n_B(S,\op_{i=1}^nQ_i)\longrightarrow\underline{\Ext}^n_B(S,S)$ be the map induced by $\pi$. Then $\pi^n_{ext}\circ\iota^n_{ext}=id$. For any element $x\in D^n=\underline{\Ext}^n_B(S,S)$, we have $y:=\iota^n_{ext}(x)\in\underline{\Ext}^n_B(S,\op_{i=1}^nQ_i)=\op_{i=1}^n\Ext_B^n(S,Q_i)$. Hence $y$ can be written as $(y_1,y_2,\dots,y_n)$ for some $y_i\in \Ext_B^n(S,B_0)=U^n$ ($i=1,\dots,n$). While $\pi\in\Hom_{B_0}(\op_{i=1}^nQ_i,S)=\underset{\text{$n$ copies}}\op S$, $\pi$ is corresponding to a sequence of elements $(s_1,\dots,s_n)$ of $S$. Now $x=\pi^n_{ext}\circ\iota^n_{ext}(x)=\pi^n_{ext}(y)=\sum_{i=1}^ns_iy_i$, where the multiplication $s_iy_i$ is the Yoneda product of the extensions, that is to say, $x$ may be written as a sum of multiplications of some elements in $V^0$ and some elements in $U^n$. Summary, we have proved that if $E$ is generated in ext-degrees 0, 1 and 2, then the graded matrix algebra
$$\left(
                                                                                  \begin{array}{cc}
                                                                                    E & U \\
                                                                                    V & D \\
                                                                                  \end{array}
                                                                                \right)
$$ is also generated in degrees 0, 1 and 2, or equivalently, the graded algebra $\underline{\Ext}^*_B(A_0,A_0)$ is generated in ext-degrees 0, 1 and 2. Hence the graded algebra $\underline{\Ext}^*_A(A_0,A_0)$ is generated in ext-degree 0, 1 and 2. Now applying \cite[Theorem 4.1]{GMMZ}, we obtain that $A$ is $N$-Koszul. \qed

At this moment, we don't know whether the converse of the theorem above is true. However, for some special cases, the converse holds.

Let $A=\op_{n\ge0}A_n$ be a positively graded Hopf algebra such that
$A_0$ is a semisimple and cosemisimple Hopf
algebra, $A_i$ is finite dimensional for all $i\ge0$ and $A_iA_j=A_{i+j}$. Let $H=A_0$.  Then the natural projection
$p:\, A\longrightarrow H$ is a Hopf algebra map; and $A$ becomes a graded right
$H$-comodule algebra through the projection $p$. Let $B=A^{coH}$. Then $B$ is a positively graded algebra
with $B_0\cong k$. Now we assume that $A/B$ is a graded $H$-Galois
extension. View $A_0$ as a right $A$-module via the projection $p$. We have the following corollary.

\begin{cor}\label{cor2} With the notation as above, $A$ is a Koszul algebra
if and only if the coinvariant subalgebra $B$ is.
\end{cor}
\proof From the proof of the theorem above we see that the graded algebra $\Ext_B^*(A_0,A_0)$ is generated
in ext-degrees 0, 1 and 2 if and only if $\Ext^*_A(A_0,A_0)$ is generated in
ext-degrees 0, 1 and 2. Since $B_0=k$, $A_0$ as a graded $B$-module is a finite sum of
copies of $k$. Hence $\Ext_B^*(A_0,A_0)$ is isomorphic to a matrix algebra of the graded algebra $\Ext_B^*(k,k)$. Hence $\Ext_B^*(A_0,A_0)$ is generated in ext-degrees 0, 1 and 2 if and only if $\Ext_B^*(k,k)$ is. Also from the proof of the theorem above, we see that $B$ is a homogeneous algebra. Now the result follows from \cite[Theorem 4.1]{GMMZ}. \qed

Let $A$ be a cocommutative Hopf algebra over an algebraic closed field of characteristic $0$. It is well-known
that $A\cong U(P(A))\#kG$ as Hopf algebras, where $P(A)$ is the Lie
algebra of primitive elements of $A$ and $U(P(A))$ is the universal
enveloping algebra of $P(A)$, and $G$ is the group of group-like
elements of $A$. Let $gr(A)$ be the graded Hopf algebra associated
to the coradical filtration of $A$. Then $gr(A)\cong
gr(U(P(A)))\#kG$. In particular, $gr(A)$ is a Galois (graded)
extension over $kG$. If $\dim(P(A))<\infty$, then $gr(U(P(A)))$ is
isomorphic to a polynomial algebra. Hence we have the following
corollary, which results in that a cocommutative  Hopf algebra is a PBW-deformations of a Koszul algebra. By the Koszul property, we may lift certain homological properties of $gr(A)$ to $A$, e.g. Calabi-Yau property, see \cite{HVZ}.

\begin{cor} Let $A$ be a cocommutative Hopf algebra over an algebraic closed field of characteristic $0$ such that $G$ is a finite group and $\dim(P(A))<\infty$. Then $gr(A)$ is a Koszul algebra.
\end{cor}

We end this note with an example of a cocommutative Hopf algebra
which satisfies Corollary \ref{cor2}.

\begin{exa} {\rm Let $Q$ be the following quiver.}
\vspace{3mm}
\unitlength=1.2pt
$$\begin{picture}(200,28)(0,0)
\put(80,10){\circle*{3}}\put(120,10){\circle*{3}}\put(75,10){\makebox(0,0){$^0$}}\put(125,10){\makebox(0,0){$^1$}}
\put(100,30){\makebox(0,0){$x_0$}}\put(100,11){\makebox{$y_0$}}\put(100,0){\makebox{$x_1$}}\put(100,-18){\makebox{$y_1$}}
\put(121,14){$\vector(4,-3){0}$}\put(120,12){$\vector(4,-3){0}$}\put(80,8){$\vector(-4,3){0}$}\put(79,6){$\vector(-4,3){0}$}
\qbezier(80,14)(100,38)(120,14)\qbezier(80,12)(100,25)(120,12)
\qbezier(80,8)(100,-15)(120,8)\qbezier(80,6)(100,-28)(120,6)
\end{picture}$$
\end{exa}
\vspace{5mm}
As an algebra, let $A=kQ/I$ where the ideal $I$ is generated by the
relations $\{x_0y_1-y_0x_1,x_1y_0-y_1x_0\}$. Denote by $e_0$ and
$e_1$ the idempotents corresponding to the vertices. We define a
coproduct and a counit on $A$ so that it becomes a Hopf algebra. It
is enough to define the comultiplication and the counit on the
generators.
$$\begin{array}{ccl}
    \Delta(e_0) & = & e_0\ot e_0+e_1\ot e_1, \\
    \Delta(e_1) & = & e_1\ot e_0+e_0\ot e_1, \\
    \Delta(x_0) & = & e_0\ot x_0+e_1\ot x_1+x_0\ot e_0+x_1\ot e_1,\\
    \Delta(y_0) & = & e_0\ot y_0+e_1\ot y_1+y_0\ot e_0+y_1\ot e_1,\\
    \Delta(x_1) & = & e_1\ot x_0+x_1\ot e_0+e_0\ot x_1+x_0\ot e_1,\\
    \Delta(y_1) & = & e_1\ot y_0+y_1\ot e_0+e_0\ot y_1+y_0\ot e_1,
  \end{array}
$$ and $$\begin{array}{l}
           \varepsilon(e_0)=1,\qquad \varepsilon(e_1)=0,\\
           \varepsilon(x_0)=\varepsilon(x_1)=\varepsilon(y_0)=\varepsilon(y_1)=0.
         \end{array}
$$
One can check that $A$ is a cocommutative bialgebra. In fact, $A$ is also a Hopf algebra, and the antipode
is given by $$\begin{array}{l}
                S(e_0)=e_0,\qquad\quad S(e_1)=e_1,\\
                S(x_0)=-x_1,\qquad S(y_0)=-y_1,\\
                S(x_1)=-x_0,\qquad S(y_1)=-y_0.
              \end{array}
$$
It is clear that $H=A_0\cong k\mathbb{Z}_2$ is a semisimple and cosemisimple Hopf algebra. One can also check that $B=A^{coH}\cong
k[x,y]$ is a Koszul algebra. Thus $A$ is a Koszul algebra.

\subsection*{Acknowledgement}  The first named author is supported by an FWO-grant and NSFC (No. 10801099).

\bibliography{}

\end{document}